\documentclass[12pt]{article}
\usepackage{amsmath,amssymb,amsthm,amsfonts,hhline,color}
\usepackage{bbm}
\usepackage[titletoc,title]{appendix}

\usepackage[dvipdfmx]{hyperref}
\usepackage{comment}
\usepackage{color}
\usepackage[dvipdfmx]{graphicx}

\newcommand{\dfr}[2]{\dfrac{#1}{#2}}
\newcommand{\cd}{\cdot}

\newcommand{\dsum}{\displaystyle \sum}

\renewcommand{\l}{\left}
\renewcommand{\r}{\right}

\newcommand{\la}{\langle}
\newcommand{\ra}{\rangle}

\newcommand{\abs}[1]{\lvert{#1}\rvert}

\newcommand{\Z}{\mathbb{Z}}
\newcommand{\C}{\mathbb{C}}
\newcommand{\R}{\mathbb{R}}

\newcommand{\F}{\mathbb{F}}
\newcommand{\M}{\mathbb{M}}

\newcommand{\End}{\mathrm{End}}

\newcommand{\aut}{\mathrm{Aut}\,}
\newcommand{\wt}{\mathrm{wt}}

\renewcommand{\hom}{\mathrm{Hom}}

\newcommand{\ad}{\mathrm{ad}\,}

\newcommand{\pii}{\pi \sqrt{-1}\, }

\newcommand{\w}{\omega}
\newcommand{\vac}{\mathbbm{1}}

\newcommand{\Fi}{\mathrm{Fi}}
\newcommand{\mymid}{\,|\,}

\makeatletter
\@addtoreset{equation}{section}

\makeatother

\theoremstyle{plain}

\newtheorem{thm}{Theorem}[section]
\newtheorem{prop}[thm]{Proposition}
\newtheorem{lem}[thm]{Lemma}

\theoremstyle{definition}
\newtheorem{df}[thm]{Definition}

\newtheorem{exam}[thm]{Example}

\newtheorem{rem}[thm]{Remark}
\newtheorem{prob}[thm]{Problem}
\newtheorem{question}[thm]{Question}

\newcommand{\sfr}[2]{\leavevmode\kern-.05em
  \raise.5ex\hbox{\the\scriptfont0 #1}\kern-.1em
  /\kern-.15em\lower.25ex\hbox{\the\scriptfont0 #2}\kern.02em}

\newcommand{\shf}{\sfr{1}{2}}
\DeclareMathOperator*{\tensor}{\otimes}
\DeclareMathOperator*{\fusion}{\boxtimes}

\newcommand{\com}{\mathrm{Com}}

\renewcommand{\ker}{\mathrm{Ker}}
\renewcommand{\o}{\mathrm{o}}

\title{3-transposition groups arising in VOA theory}

\author{Hiroshi Yamauchi\footnote{Partially supported by Grant-in-Aid for Scientific Research (C) 19K03409.}
  \medskip\\
  {\small \it Department of Mathematics,
  Tokyo Woman's Christian University}\\
  {\small \it 2-6-1 Zempukuji, Suginami-ku, Tokyo 167-8585, Japan}\\
  {\small e-mail: \texttt{yamauchi@lab.twcu.ac.jp}}
}
\date{}

\begin{document}

\maketitle

\begin{abstract}
We review 3-transposition groups arising in vertex operator algebra theory.
One can construct a commutative algebra called the Matsuo algebra out of a 3-transposition group.
Some 3-transposition groups arise as automorphism groups of vertex operator algebras 
via Matsuo algebras but there exist some 3-transposition groups 
which do not arise through Matsuo algebras.
We will exhibit examples of those groups together with VOAs.
This article is based on the author's talk at the seminar on 
``Majorana, Axial, Vertex Algebras and the Monster (MAVAM)'' held online on June 4th, 2021.
\end{abstract}

\pagestyle{plain}

\baselineskip 6mm

\section{Introduction}
An \emph{$n$-transposition group} is a pair ($G$, $I$) of a group $G$ and its subset 
$I$ of involutions such that $I$ is a union of conjugacy classes of $G$ 
and the product of any two involutions in $I$ has order at most $n$.
3-transposition groups are first considered by Fischer \cite{Fi} and three sporadic groups 
called Fischer groups are discovered on the way of classification of finite 3-transposition groups.
Cuypers and Hall \cite{CH} have investigated the classification program further.

Given a 3-transposition group ($G$,$I$), one can define a graph structure on $I$ by adjoining 
two distinct involutions in $I$ if and only if their product has order 3.
The permutation representation of $G$ on the graph $I$ has rank 3 under the assumptions that 
$Z(G)=O_2(G)=O_3(G)=1$ and $G'=G''$ (cf.~\cite{As,Fi,CH}).
The study of center-free 3-transposition groups is actually equivalent to that of the 
graph structures defined as above.
As a linearization of the graph structure, one can define a commutative algebra called 
the Matsuo algebra associated with a 3-transposition group which is linearly spanned 
by idempotents indexed by transpositions (cf.~Definition \ref{df:3.1}).
Therefore, if we realize the Matsuo algebra associated with a 3-transposition group as 
a substructure of some algebraic structures, then we may obtain a natural action of the 
3-transposition group on the algebraic  structures under consideration.
This is a way how 3-transposition groups arise in the theory of vertex operator algebras.

A vertex operator algebra (VOA) is an infinite dimensional vector space equipped with 
infinitely many bilinear operations indexed by integers satisfying a certain set of axioms 
(cf.~Definition \ref{df:4.1}).
Given a VOA $V$ and its sub VOA $U$, one can decompose $V$ as a $U$-module and if 
$U$ is rational then $V$ is a direct sum of irreducible $U$-submodules.
In this case fusion rules among $U$-modules constrain binary operations in $V$ 
so that if the fusion rules among irreducible $U$-modules have a cyclic symmetry then 
this symmetry gives rise to define an automorphism of $V$ through its decomposition as 
a $U$-module.
This is a typical way to obtain automorphisms of VOAs and in practice $U$ is often chosen 
to be a subalgebra with a small number of generators such as a Virasoro VOA or a $W_3$-algebra 
since they are easy to handle.
If $U$ is a simple Virasoro VOA in the BPZ series \cite{BPZ} then its fusion rules 
have a $\Z_2$-symmetry and the involutive automorphism defined by $U$ is called 
the \emph{Miyamoto involution} (cf.~\cite{M1}).
The most fundamental example is the Virasoro VOA with central charge 1/2.
A Virasoro vector is called an \emph{Ising vector} if it generates 
a simple Virasoro sub VOA with central charge 1/2.
The sub VOA generated by two Ising vectors are studied in \cite{S} and it is shown that 
the product of Miyamoto involutions associated with two Ising vectors is bounded by 6 
which is nowadays known as Sakuma's theorem.
Thus in general we cannot directly obtain 3-transposition groups if we consider 
automorphism groups generated by Miyamoto involutions associated with Ising vectors.

This article is an exposition of 3-transposition groups arising in VOA theory.
We will mainly consider VOAs of OZ-type over $\R$ with positive definite 
invariant bilinear forms.
The degree 2 subspace of a VOA of OZ-type carries a structure of a commutative but 
in general non-associative algebra called the Griess algebra with symmetric invariant bilinear 
form (cf.~\cite{FLM,G}).
Let $V$ be a VOA of OZ-type over $\R$.
It is known that if $x\in V_2$ is an idempotent of the Griess algebra of $V$ 
then $2x$ is a Virasoro vector of $V$ with central charge $8(x\mymid x)$ (cf.~\cite{M1}).
Sakuma's theorem \cite{S} mentioned above describes subalgebras of the Griess algebra 
generated by two Ising vectors which correspond to idempotents with squared norm 1/16.
The notions of Majorana algebras \cite{Iv} and axial algebras \cite{HRS} are introduced 
by axiomatizing the role of Ising vectors in the Griess algebra.
Recall that a Matsuo algebra is a commutative algebra linearly spanned by idempotents of 
fixed squared norm.
Those idempotents are called axial vectors.
If the Griess algebra $V_2$ contains a subalgebra isomorphic to (the non-degenerate quotient of) 
a Matsuo algebra and each axial vector generates a rational Virasoro sub VOA 
then Miyamoto involutions associated with axial vectors generate a 3-transposition group.
In this article we will review some families of Matsuo algebras which arise as Griess algebras
of VOAs and related topics like the classification problems and the uniqueness problems.
On the other hand, there are constructions of 3-transposition groups as automorphism groups 
of VOAs based on Ising vectors which do not arise from Matsuo algebras.
We will also review such examples.

The organization of this article is as follows.
In Section 2 we will review 3-transposition groups based on \cite{CH,K}.
In Section 3 we will review Matsuo algebras based on \cite{Ma}.
Note that the published version of \cite{Ma} and the one posted to the arXiv are different.
For \cite{Ma} we refer to the arXiv version in this article.
In Section 4 we will review the notion of vertex operator algebras and Miyamoto involutions.
Section 5 is the main body of this article and we will review several constructions of 
3-transposition groups in VOA theory except for Section 5.5 where a construction of 
4-transposition groups will be explained.
We will close this article by proposing some open questions in Section 6.

\paragraph{Acknowledgement.}
This article is based on the author's talk at the seminar on 
``Majorana, Axial, Vertex Algebras and the Monster (MAVAM)'' held online on June 4th, 2021.
The author wishes to thank the organizers of MAVAM seminar for the opportunity to give 
a talk there.

\section{3-transposition groups}

Recall the notion of 3-transposition groups (cf.~\cite{As,CH,Fi,K}).

\begin{df}\label{df:2.1}
A \emph{3-transposition group} is a pair ($G$, $I$) of a group $G$ and its subset 
$I$ of involutions called \emph{transpositions} such that $I$ is closed under the conjugation 
and the product of any two involutions in $I$ has order at most 3.
A subgroup of $G$ generated by a subset of $I$ is called an \emph{$I$-subgroup} which is again 
a 3-transposition group.
\end{df}

In the definition above, some authors indeed require $I$ to be a single conjugacy class of $G$.
To include some subgroups such as $\mathfrak{S}_k\times \mathfrak{S}_{n-k}$ of $\mathfrak{S}_n$, 
we admit $I$ to be a union of conjugacy classes.

\begin{exam}
Here are some typical examples of 3-transposition groups.
\\
(1) The symmetric group $G=\mathfrak{S}_n$ with $I=\{ (i~j) \mid 1\leq i<j\leq n\}$.
\\
(2) Weyl groups of ADE types.
\\
(3) Orthogonal groups over $\F_2$ where the 3-transpositions are the transvections.
\\
(4) Symplectic groups over $\F_2$ where the 3-transpositions are the transvections.
\\
(5) Orthogonal groups over $\F_3$ where the 3-transpositions are the reflections.
\\
(6) Unitary groups over $\F_4$ where the 3-transpositions are the transvections.
\end{exam}

We have the following isomorphisms between several 3-transposition groups 
with different base fields:
\[
\begin{array}{l}
  U_4(2)\cong O_5(3) < O_5(3){:}2\cong O^-_6(2){:}2
  \cong W(E_6),~~ 
  S_6(2)\cong W(E_7)/\{ \pm 1\}, ~~
  \medskip\\
  O_8^+(2)\cong W(E_8)/\{\pm 1\},~~~
  O_6^+(2){:}2 \cong \mathfrak{S}_8,~~~
  O_4^-(3){:}2\cong S_4(2)\cong \mathfrak{S}_6,~~~
  O_4^-(2){:}2 \cong \mathfrak{S}_5.
\end{array}
\]
Here we use the notation as in \cite{ATLAS} for simple groups 
and $W(X)$ denotes the Weyl group of type $X$.

The study of 3-transposition groups was initiated by Fischer and under certain assumptions 
he classified finite 3-transposition groups up to centers as follows.

\begin{thm}[\cite{Fi, CH}]\label{thm:2.3}
Let $G$ be a finite 3-transposition group generated by a single conjugacy class of 
3-transpositions.
If $O_2(G)O_3(G)<Z(G)$ then up to the center, $G$ is isomorphic to one of the following:
\\
\textup{(1)} Symmetric groups $\mathfrak{S}_n$, $n\geq 5$.
\\
\textup{(2)} Symplectic groups $S_{2n}(2)$, $n\geq 2$.
\\
\textup{(3)} Orthogonal groups $O^\varepsilon_{2n}(2){:}2$ 
where $n\geq 3$ if $\varepsilon=+$ and $n\geq 2$ if $\varepsilon=-$.
\\
\textup{(4)} Orthogonal groups $O^\varepsilon_{2n}(3){:}2$
where $n\geq 3$ if $\varepsilon=+$ and $n\geq 2$ if $\varepsilon=-$.
\\
\textup{(5)} Orthogonal groups $O^\pm_{2n+1}(3)$ or $O^\pm_{2n+1}(3){:}2$ with $n\geq 2$ and 
$\varepsilon=\pm$.
\\
\textup{(6)} Unitary groups $U_n(2)$, $n\geq 4$.
\\
\textup{(7)} One of 5 exceptionals: $\Fi_{22}$, $\Fi_{23}$, $\Fi_{24}$, 
$O_8^+(2){:}\mathfrak{S}_3$ or $O_8^+(3){:}\mathfrak{S}_3$.
\end{thm}

Cupyers and Hall classified the structure $G/Z(G)$ of a 3-transposition group $G$ without 
the assumptions that $G$ is finite and $O_2(G)O_3(G)<Z(G)$ in \cite{CH}.
In their classification, 3-transposition groups are divided into four classes as follows. 
Set 
\begin{equation}\label{eq:2.1}
  H:=\la a,b,c \mid a^2=b^2=c^2=(ab)^3=(bc)^3=(ca)^3=(abac)^3=1\ra .
\end{equation}
Then $Z(H)=\la (abc)^2\ra\cong 3$ and $H\cong 3^{1+2}{:}2$.
The group $H$ is denoted by $\mathrm{SU}_3(2)'$ in \cite{CH}.
A 3-transposition group $(G,I)$ is called of \emph{symplectic type} if $G$ does not contain 
an $I$-subgroup isomorphic to $H$, and $G$ is of \emph{orthogonal type} if it contains an 
$I$-subgroup isomorphic to $H$ but no $I$-subgroup isomorphic to $2^{1+6}{:}H$, 
and $G$ is of \emph{unitary type} if it contains an $I$-subgroup isomorphic to $2^{1+6}{:}H$ 
but no $\mathrm{P}\Omega_8^+(2){:}\mathfrak{S}_3$, and $G$ is of \emph{sporadic type} 
if $G$ does contain an $I$-subgroup isomorphic to 
$\mathrm{P}\Omega_8^+(2){:}\mathfrak{S}_3$ (cf.~[CH]).
Among the groups in Theorem \ref{thm:2.3}, those in (1), (2) and (3) are of symplectic type, 
those in (4) and (5) are of orthogonal type, those in (6) are of unitary type and 
those in (7) are of sporadic type. 
Note that $O^{\mp}_{2n}(2)<S_{2n}(2)<O^{\pm}_{2(n+1)}(2)$,  
$\mathfrak{S}_{2n+1}<S_{2n}(2)$ and $\mathfrak{S}_{2n+2}<S_{2n}(2)$ 
so that the groups in (1), (2), (3) of Theorem \ref{thm:2.3} are called of symplectic type.

\begin{rem}
In the definition of $H$ in \eqref{eq:2.1}, if we replace the relation $(abac)^3=1$ by 
$abac=1$ then we obtain $\mathfrak{S}_3$, and if we replace the relation $(abac)^3=1$ by 
$(abac)^2=1$ then we obtain $\mathfrak{S}_4$.
\end{rem}

\section{Matsuo algebras}

Let $(G,I)$ be a 3-transposition group.
We define a binary relation on $I$ such that $a\sim b$ for $a$, $b\in I$ 
if and only if $ab$ has order 3, that is, we have $a^b=b^a$ in $I$, 
and in this case we set $a\circ b:=a^b=b^a\in I$.
This binary relation amounts to define an irreflexive undirected graph structure on $I$.
We say $(G,I)$ (or just $G$) is \emph{connected} or \emph{indecomposable} 
if $I$ is connected with respect to the adjacency relation defined by $a\sim b$.
This is equivalent to say that $I$ is a single conjugacy class of involutions.

Based on the graph structure on $I$, we define the Matsuo algebra associated with 
a 3-transposition group as follows.

\begin{df}[\cite{Ma}]\label{df:3.1}
  For $\alpha$, $\beta\in \R$, set $B_{\alpha,\beta}(G):=\oplus_{i\in I} \R x^i$ 
  where $\{ x^i\mid i\in I\}$ is a formal basis and define the product and 
  the bilinear form as follows.
  \begin{equation}\label{eq:3.1}
    x^i x^j :=
    \begin{cases}
      2x^i  & \mbox{if $i=j$,}
      \medskip\\
      \dfr{\alpha}{2} (x^i+x^j-x^{i\circ j}) & \mbox{if $i\sim j$,}
      \medskip\\
      0 & \mbox{otherwise,}
    \end{cases}
    ~~~~~~~~~
    (x^i \mymid x^j):=
    \begin{cases}
      \dfr{\beta}{2} & \mbox{if $i=j$,}
      \medskip\\
      \dfr{\alpha\beta}{8} & \mbox{if $i\sim j$,}
      \medskip\\
      0 & \mbox{otherwise.}
    \end{cases}
  \end{equation}
  Then $B_{\alpha,\beta}(G)$ becomes a commutative algebra with symmetric invariant bilinear form.
  We call $B_{\alpha,\beta}(G)$ the \emph{Matsuo algebra} associated with $(G,I)$, and we also 
  call each $x^i$ with $i\in I$ an \emph{axial vector}.
  By the invariance, the radical of the bilinear form becomes an ideal of $B_{\alpha,\beta}(G)$. 
  We call the quotient algebra of $B_{\alpha,\beta}(G)$ by the radical of the bilinear form 
  the \emph{non-degenerate quotient}. 
\end{df}

\begin{rem}
By definition, $x^i/2$ is an idempotent of $B_{\alpha,\beta}(G)$ so that a Matsuo algebra is 
linearly spanned by idempotents.
\end{rem}

\paragraph{Unity.}
Suppose the set $I$ of 3-transpositions is not a single conjugacy class, i.e., 
$I=I_1\sqcup I_2$ where $I_1$ and $I_2$ are disconnected with respect to the adjacency relation.
Then $G$ is isomorphic to a central product $G_1*G_2$ of $G_1=\la I_1\ra$ and $G_2=\la I_2\ra$. 
Correspondingly, we have a decomposition into a direct sum 
\[
  B_{\alpha,\beta}(G)=B_{\alpha,\beta}(G_1)\oplus B_{\alpha,\beta}(G_2),
\]
of ideals which is also an orthogonal sum.
Therefore, in the study of the Matsuo algebra $B_{\alpha,\beta}(G)$, we may assume 
that $G$ is indecomposable, that is, $G$ is generated by a single conjugacy class of 
3-transpositions.
In particular, $B_{\alpha,\beta}(G)$ is unital if and only if the Matsuo algebra associated with 
each indecomposable component of $G$ has the unity.

Suppose $G$ is indecomposable. 
Then $k:=\#\{ j\in I \mid i\sim j\}$ is independent of $i\in I$ and the vector 
\begin{equation}\label{eq:3.2}
  \w:=\dfr{4}{k\alpha+4}\dsum_{i\in I} x^i
\end{equation}
satisfies $\w x^i=2x^i$ and $(\w \mymid x^i)=(x^i\mymid x^i)=\beta/2$ for all $i\in I$ 
so that $\w /2$ gives the unity of $B_{\alpha,\beta}(G)$ if $k\alpha+4\ne 0$.

\begin{rem}
When the Matsuo algebra $B_{\alpha,\beta}(G)$ (or its quotient) is realized by a VOA 
as its Griess algebra, then twice the unity of $B_{\alpha,\beta}(G)$ corresponds to
the conformal vector of the VOA. 
In this case, each axis $x^i$ gives a Virasoro vector with central charge $\beta$.
\end{rem}

\paragraph{Involutions.}
It is obvious that $G$ acts on $B_{\alpha,\beta}(G)$ by conjugation.
Namely, there exists a unique group homomorphism
\begin{equation}\label{eq:3.3}
\begin{array}{cccc}
  \sigma : &G & \longrightarrow & \aut B_{\alpha,\beta}(G)
  \medskip\\
  & a & \longmapsto & \sigma_a : x^i\mapsto x^{aia^{-1}}
\end{array}
\end{equation}
such that $\sigma_i(x^j)=x^{i\circ j}$ if $i\sim j$ and $\sigma_i(x^j)=x^j$ otherwise
for $i$, $j\in I$.
Clearly $\ker\,\sigma= Z(G)$ and $\sigma$ is injective if $G$ is center-free.

Suppose $\alpha \ne 2$.
Then the adjoint action $\ad x^i$ of $x^i$ is semisimple and has three eigenvalues $0$, $2$ 
and $\alpha$ on $B_{\alpha,\beta}(G)$.
For, clearly $x^i$ is an eigenvector of $\ad x^i$ with eigenvalue 2.
Let $j\in I$ be distinct from $i$.
If $i\not\sim j$, then $(\ad x^i)\, x^j=x^ix^j=0$ so that $x^j$ belongs to the eigenspace of 
eigenvalue 0.
If $i\sim j$, then $i\circ (i\circ j)=j$ and one has 
\[
  x^i(x^j+x^{i\circ j})=\alpha x^i,~~~
  x^i(x^j-x^{i\circ j})=\alpha(x^j-x^{i\circ j}), 
\]
from which it follows that $\ad x^i$ is semisimple and has eigenvalues 0, 2 and $\alpha$.
One can also verify that $\ker\, (\ad x^i -2)=\R x^i$ and $\sigma_i$ acts by the identity on 
$\ker\, \ad x^i$ and by $-1$ on $\ker\, (\ad x^i-\alpha)$.
Therefore, we have a decomposition of $B_{\alpha,\beta}(G)$ such that
\begin{equation}\label{eq:3.4}
\begin{array}{ccccccc}
  B_{\alpha,\beta}(G) &= & \R x^i & \oplus & \ker\, (\ad x^i) & \oplus & \ker\, (\ad x^i-\alpha),
  \smallskip\\
  \ad x^i & : & 2 && 0 && \alpha,
  \smallskip\\
  \sigma_i & : & 1 && 1 && -1.
\end{array}
\end{equation}
If we set $B_{\alpha,\beta}^{i,\pm}(G):= \{ v\in B_{\alpha,\beta}(G) \mid \sigma_i(v)=\pm v\}$
then we have $B_{\alpha,\beta}(G)=B_{\alpha,\beta}^{i,+}(G)\oplus B_{\alpha,\beta}^{i,-}(G)$
with 
\[
  B_{\alpha,\beta}^{i,+}(G) = \R x^i  \oplus  \ker\, (\ad x^i),~~~~
  B_{\alpha,\beta}^{i,-}(G)=\ker\, (\ad x^i-\alpha).
\]
Therefore, the action of $G$ on $B_{\alpha,\beta}(G)$ can be described 
by the adjoint action of $B_{\alpha,\beta}(G)$.
The involution $\sigma_i$ defined as in \eqref{eq:3.4} is called the 
\emph{Miyamoto involution} associated with the axial vector $x^i$ of $B_{\alpha,\beta}(G)$.

\section{Vertex operator algebras}

Recall the notion of vertex operator algebras (VOAs).

\begin{df}\label{df:4.1}
A vertex operator algebra is a $\Z_{\geq 0}$-graded vector space $V=\bigoplus_{n\geq 0} V_n$ 
over $\C$ equipped with bilinear products $V\tensor V\ni a\tensor b\mapsto a_{(n)}b\in V$ 
indexed by $n\in \Z$ satisfying the following axioms.
\begin{enumerate}
\item 
For any $a$, $b\in V$, there exists $N\in \Z$ such that $a_{(i)}b=0$ for $i\geq N$.

\item
There exists a special element $\vac \in V_0$ called the \emph{vacuum vector} of $V$ 
such that for any $a\in V$ one has $a_{(n)}\vac=\delta_{n,-1}a$ for $n\geq -1$.

\item
There exists a special element $\w \in V_2$ called the \emph{conformal vector} of $V$ 
such that if we set $L(n)=\w_{(n+1)}\in \End(V)$ for $n\in \Z$ then we have 
\begin{equation}\label{eq:4.1}
  [L(m),L(n)]=(m-n)L(m+n)+\delta_{m+n,0}\dfr{m^3-m}{12}c_V, 
\end{equation}
where $c_V\in \C$ is called the \emph{central charge} of $V$. 
We also have 
\[
  [L(-1),a_{(m)}]=-m\,a_{(m-1)}
\] 
for any $a\in V$ and $m\in \Z$ and $V_n= \ker\, (L(0)-n)$ with $\dim V_n<\infty$.

\item
The following Borcherds identity holds for any $a$, $b$, $c\in V$ and $k$, $m$, $n\in \Z$.
\[
  \dsum_{i=0}^\infty \binom{m}{i}\l( a_{(k+i)}b\r)_{(m+n-i)}c
  =\dsum_{i=0}^\infty (-1)^i \binom{k}{i}\!\l\{ a_{(m+k-i)}b_{(n+i)} -(-1)^k b_{(n+k-i)}a_{(m+i)}\r\} c.
\]
\end{enumerate}
\end{df}

There are various equivalent but apparently different formulations of 
the vertex operator algebra.
The set of axioms above is an apparently minimum one (cf.~\cite{MaN}).
We refer the readers to \cite{FHL,FLM,Ka,LeLi,MaN} for more details.

If $a \in V_n$ is homogeneous with respect to the grading of $V$ then we define $\wt(a):=n$ 
and call it the \emph{weight} of $a$.
It follows from the axioms above that for homogeneous $a$, $b\in V$, 
the product $a_{(n)}b$ is also homogeneous with weight $\wt(a)+\wt(b)-n-1$ for $n\in \Z$.
In particular, the \emph{zero-mode operator} $\o(a)=a_{(\wt(a)-1)}\in \End (V)$ preserves 
each graded subspace of $V$.

\paragraph{Griess algebra.}

In this article we will mainly consider VOAs of OZ-type defined over $\R$.
A VOA $V=\bigoplus_{n\geq 0} V_n$ over $\R$ is called of \emph{OZ-type}\footnote{%
OZ stands for ``One-Zero''. This notion is introduced by Bob Griess.}
 if $V_0=\R \vac$ and $V_1=0$.
In this case we can define a product and a bilinear form on the degree 2 subspace $V_2$ 
as follows.
Let $a$, $b\in V_2$.
We define
\begin{equation}\label{eq:4.2}
  ab:=\o(a)b=a_{(1)}b \in V_2,~~~~~~~~
  (a\,|\,b)\vac := a_{(3)}b \in V_0=\R \vac.
\end{equation}
Then $V_2$ forms a commutative algebra with symmetric invariant bilinear form 
(cf.~Proposition 8.2.1 of \cite{MaN}, see also Section 8.9 of \cite{FLM}).
The algebra $V_2$ is called the \emph{Griess algebra} of $V$.

\begin{rem}\label{rem:4.2}
A VOA of OZ-type has a unique symmetric invariant bilinear form defined by 
\[
  (\vac\,|\vac):=1 
  ~~\mbox{ and } ~~
  (a\, |\, b):=\dsum_{i\geq 0} \dfr{(-1)^{\wt(a)}}{i!}(\vac \, |\, 
  \underbrace{(L(1)^ia)_{(2\wt(a)-i-1)}b}_{\in V_0=\R \vac}\,).
\]
For details, see \cite{FHL,Li}.
The bilinear form on $V_2$ defined in \eqref{eq:4.2} is a restriction of the form above 
to the weight 2 subspace of $V$.
\end{rem}

In the rest of the article, we always assume that \textbf{our VOA $V$ is defined over $\R$, 
of OZ-type, and the invariant bilinear form on $V$ is positive definite} unless otherwise stated.
Due to the positivity, it follows that every sub VOA of $V$ is always simple.

\medskip

A vector $a\in V_2$ is called a \emph{Virasoro vector} with central charge $c_a\in \R$ 
if the modes $L^a(n)=a_{(n+1)}\in \End(V)$, $n\in \Z$, satisfy the commutation relation
\begin{equation}\label{eq:4.3}
  [L^a(m),L^a(n)]=(m-n)L^a(m+n)+\delta_{m+n,0}\dfr{m^3-m}{12}c_a
\end{equation}
of the Virasoro algebra.
By \eqref{eq:4.1}, the conformal vector of $V$ is a Virasoro vector 
(but not vice versa in general).
For a VOA of OZ-type, Virasoro vectors corresponds to idempotents of the Griess algebra.

\begin{lem}[\cite{M1}]
  A vector $a\in V_2$ is an idempotent of the Griess algebra of $V$ if and only if 
  $2a$ is a Virasoro vector with central charge $8(a\mymid a)$.
\end{lem}


A sub VOA $\la a\ra$ generated by a Virasoro vector $a\in V_2$ is called a Virasoro VOA.
Thanks to the assumption of the positivity, a Virasoro sub VOA is always simple in our setting.

\paragraph{Unitary series.}

Recall the unitary series of the Virasoro algebra.
\begin{equation}\label{eq:4.4}
\begin{array}{l}
  c_m:=1-\dfr{6}{(m+2)(m+3)}~~~~(m=1,2,3,\dots)
  \medskip\\
  h_{r,s}^{(m)}:=\dfr{(r(m+3)-s(m+2))^2-1}{4(m+2)(m+3)},~~~
  1\leq r\leq m+1,~~ 1\leq s \leq m+2
\end{array}
\end{equation}
In this article $L(c_m,h_{r,s}^{(m)})$ denotes the irreducible highest weight modules over 
the Virasoro algebra with central charge $c_m$ and highest weight $h_{r,s}^{(m)}$.
It was shown in \cite{FZ} that $L(c_m,0)$ forms a simple and rational Virasoro VOA and 
its irreducible $L(c_m,0)$-modules are provided by highest weight representations 
$L(c_m,h_{r,s}^{(m)})$, $1\leq s\leq r\leq m+1$.
Note that we have redundancy $h_{r,s}^{(m)}=h_{m+2-r,m+3-s}^{(m)}$.
The fusion rules among irreducible $L(c_m,0)$-modules are computed in \cite{W} as follows.
\begin{equation}\label{eq:4.5}
  L(c_m,h_{r,s}^{(m)})\fusion_{L(c_m,0)} L(c_m,h_{r',s'}^{(m)})
  = \bigoplus_{1\leq i\leq I\atop 1\leq j\leq J} L(c_m,h_{\abs{r-r'}+2i-1,\abs{s-s'}+2j-1}), 
\end{equation}
where $I=\min\{ r,r',m+2-r,m+2-r'\}$ and $J=\min\{ s,s',m+3-s,m+3-s'\}$.

\paragraph{Miyamoto involutions.}

Let $x\in V$ be a Virasoro vector with central charge $c_m$.
Since $\la x\ra \cong L(c_m,0)$ is rational, we have the isotypical decomposition 
\begin{equation}\label{eq:4.6}
  V=\bigoplus_{1\leq s\leq r\leq m+1} V[x;h^{(m)}_{r,s}],~~~
  V[x;h^{(m)}_{r,s}]\cong L(c_m,h_{r,s}^{(m)})\tensor \hom_{\la x\ra}(L(c_m,h_{r,s}^{(m)}),V).
\end{equation}
Based on the decomposition above, we can define a linear automorphism $\tau_x$ of $V$ as follows.
\begin{equation}\label{eq:4.7}
  \tau_x:=(-1)^{4(m+2)(m+3)\o(x)}
  =\begin{cases} (-1)^{r+1} \mbox{~ on ~$V[x;h_{r,s}^{(m)}]$} & \mbox{if $m$ is even,}
  \medskip\\ 
  (-1)^{s+1} \mbox{~ on ~$V[x;h_{r,s}^{(m)}]$} & \mbox{if $m$ is odd.} 
  \end{cases}
\end{equation}
Then it follows from the fusion rules in \eqref{eq:4.5} that $\tau_x$ preserve the VOA structure. 

\begin{thm}[\cite{M1}]
  $\tau_x\in \aut V$.
\end{thm}

For instance, we consider the case $m=1$ in the unitary series \eqref{eq:4.4}.
Suppose $x\in V$ is a Virasoro vector with central charge $c_1=1/2$.
Then $\la x\ra \cong L(\shf,0)$ have three irreducible representations 
$L(\shf,0)$, $L(\shf,\shf)$ and $L(\shf,\sfr{1}{16})$ so that we have the following 
isotypical decomposition as in \eqref{eq:4.6}.
\begin{equation}\label{eq:4.8}
\begin{array}{ccccccc}
  V&=& V[x;0]&\oplus &V[x;\shf]&\oplus &V[x;\sfr{1}{16}]
  \medskip\\
  \tau_x &:& 1 && 1&& -1
\end{array}
\end{equation}
Since $V_2[x;0]=\R x\oplus \ker_{V_2} \o(x)$ and $V_2[x;h]=\ker_{V_2} (\o(x)-h)$ for $h=1/2$ and 
$1/16$, the corresponding decomposition of the Griess algebra looks as follows.
\begin{equation}\label{eq:4.9}
\begin{array}{ccccccccc}
  V_2&=& \R x &\oplus & \ker_{V_2} \o(x) &\oplus & \ker_{V_2} (\o(x)-1/2) & \oplus & 
  \ker_{V_2} (\o(x)-1/16)
  \medskip\\
  \o(x) &:& 2 && 0 && 1/2 && 1/16
  \medskip\\
  \tau_x &:& 1 && 1 && 1&& -1
\end{array}
\end{equation}
By \eqref{eq:4.8}, we see that the fixed point subalgebra $V^{\la \tau_x\ra}$ has a decomposition 
as follows.
\begin{equation}\label{eq:4.10}
\begin{array}{ccccc}
  V^{\la \tau_x\ra}&=& V[x;0]&\oplus &V[x;\shf]
  \smallskip\\
  \sigma_x &:& 1 && -1
\end{array}
\end{equation}
If we define the linear automorphism $\sigma_x$ as above, then again by the fusion rules 
in \eqref{eq:4.5} it follows that $\sigma_x\in \aut V^{\la \tau_x\ra}$ (cf.~\cite{M1}).
By \eqref{eq:4.9}, the action of $\sigma_x$ on the Griess algebra looks as follows.
\begin{equation}\label{eq:4.11}
\begin{array}{ccccccc}
  V_2&=& \R x &\oplus & \ker_{V_2} \o(x) &\oplus & \ker_{V_2} (\o(x)-1/2) 
  \medskip\\
  \o(x) &:& 2 && 0 && 1/2
  \medskip\\
  \sigma_x &:& 1 && 1 && -1
\end{array}
\end{equation}
We can generalize the above construction of the second automorphism as follows.

Let $x$ be a $c=c_m$ Virasoro vector of $V$ as before.
Set 
\begin{equation}\label{eq:4.12}
  P_m:=\begin{cases}~ \{ h_{1,s}^{(m)} \mid 1\leq s\leq m+2\} & \mbox{if $m$ is even,}
  \medskip\\
  ~ \{ h_{r,1}^{(m)} \mid 1\leq r\leq m+1\} & \mbox{if $m$ is odd,}
  \end{cases}
\end{equation}
and define $V[x;P_m]:=\bigoplus_{h\in P_m} V[x;h]$.
Then it follows from the fusion rules of $L(c_m,0)$-modules in \eqref{eq:4.5} that 
the subspace $V[x;P_m]$ forms a subalgebra of the fixed point subalgebra $V^{\la \tau_x\ra}$.
We define a linear automorphism $\sigma_x$ of $V[x;P_m]$ as follows.
\begin{equation}\label{eq:4.13}
  \sigma_x
  :=
  \begin{cases} 
  (-1)^{s+1} \mbox{~ on ~$V[x;h_{1,s}^{(m)}]$} & \mbox{if $m$ is even,} 
  \smallskip\\ 
  (-1)^{r+1} \mbox{~ on ~$V[x;h_{r,1}^{(m)}]$} & \mbox{if $m$ is odd.} 
  \end{cases}
\end{equation}

\begin{thm}[\cite{M1}]
  $\sigma_x \in \aut V[x;P_m]$.
\end{thm}

The involutions $\tau_x\in \aut V$ and $\sigma_x\in \aut V[x;P_m]$ are called 
\emph{Miyamoto involutions}.
Note that $\tau_x$ is always well-defined on the whole space $V$ whereas $\sigma_x$ is 
locally defined on the subspace $V[x;P_m]$.
We call a $c=c_m$ Virasoro vector $x\in V$ of \emph{$\sigma$-type on} $V$ if $V=V[x;P_m]$.
Note that if $x$ is of $\sigma$-type then $V=V^{\la \tau_x\ra}$ so that $\tau_x$ 
is trivial on $V$ and $\sigma_x$ is globally defined on $V$.

\paragraph{Sakuma's theorem.}

A $c=1/2$ Virasoro vector $x \in V$ is called an \emph{Ising vector} of $V$.
The structures of a subalgebra generated by two Ising vectors as well as the dihedral groups 
generated by associated Miyamoto involutions are almost classified by Sakuma.

\begin{thm}[\cite{S}]\label{thm:4.6}
Let $e$ and $f$ be Ising vectors of $V$.
Then the Griess algebra of the subalgebra $\la e,f\ra$ has 9 possible structures 
and it follows that $\abs{\,\tau_{e}\tau_{f}\,}\leq 6$ on $V$.
\renewcommand{\arraystretch}{1.5}
\[
\begin{array}{|c||c|c|c|c|c|c|c|c|c|}
\hline
\mbox{\textup{Type of} } \la e,f\ra & ~1\mathrm{A}~ & ~2\mathrm{A}~ & ~3\mathrm{A}~ & ~4\mathrm{A}~ & ~5\mathrm{A}~ 
& ~6\mathrm{A}~ & ~4\mathrm{B}~ & ~2\mathrm{B}~ & ~3\mathrm{C}~
\\ \hline
\textup{max}\, \abs{\tau_e\tau_f} & 1 & 2 & 3 & 4 & 5 & 6 & 4 & 2 & 3 
\\ \hline
~2^{10}(e\,|\,f)~ & 2^8 & 2^5 & 13 & 2^3 & 6 & 5 & 2^2 & 0 & 2^2
\\ \hline
\dim \la e,f\ra_2 & 1 & 3 & 4 & 5 & 6 & 8 & 5 & 2 & 3
\\ \hline
\#\mbox{\textup{ of Ising vectors}} & 1 & 3 & 3 & 4 & 5 & 7 & 5 & 2 & 3
\\ \hline
\mbox{\textup{Miyamoto type}} & \sigma & \sigma & \tau & \tau & \tau & \tau & \tau & \sigma & \tau 
\\ \hline
\end{array}
\]
\renewcommand{\arraystretch}{1}
In the table above, $\mbox{max}\,\abs{\tau_e\tau_f}$ denotes the possible maximum order of 
$\tau_e\tau_f$ on $V$ and Miyamoto type denotes the type of Ising vectors $e$ and $f$ on 
the subalgebra $\la e,f\ra$.
\end{thm}

If $\la e,f\ra$ is the subalgebra of type $nX$ in the first row of the table above, then 
we will call $\la e,f\ra$, which is also denoted by $U_{nX}$, 
the \emph{dihedral subalgebra} of type $nX$ or shortly just the \emph{$nX$-algebra}.
We refer the readers to \cite{LYY} for more information about dihedral subalgebras.
In this article, we mainly treat dihedral subalgebras of types 2A and 3A.
Here we list some properties of the 2A and 3A-algebras.

\begin{lem}\label{lem:4.7}
Let $e$ and $f$ be Ising vectors of a VOA $V$.
\\
\textup{(1)}~ $\la e,f\ra$ is a 2A-algebra if and only if  $(e\,|\,f)=2^{-5}$.
\\
\textup{(2)}~ $\la e,f\ra$ is a 3A-algebra if and only if $(e\,|\,f)=13\cd 2^{-10}$.
\end{lem}

\begin{thm}[\cite{M1}]\label{thm:4.8}
Let $\la e,f\ra\cong U_{\mathrm{2A}}$ be a 2A-subalgebra of a VOA $V$.
\\
\textup{(1)}~ $U_{\mathrm{2A}}$ has three Ising vectors  $e$, $f$ and $\sigma_{e}f=\sigma_{f}e$, 
all of them are of $\sigma$-type on $U_{\mathrm{2A}}$.
\\
\textup{(2)}~ $\tau_e\tau_f=\tau_{\sigma_e f}$ on $V$.
\\
\textup{(3)}~ $\la \tau_e,\tau_f\ra$ is an elementary abelian 2-group  of order at most 4 
(and possibly trivial).
\\
\textup{(4)}~ $\abs{\sigma_e\sigma_f}=3$ on $\la e,f\ra$ and 
$\aut U_{\mathrm{2A}}=\la \sigma_e,\sigma_f\ra \cong \mathfrak{S}_3$.
\end{thm}

\begin{thm}[\cite{M3,SY}]\label{thm:4.9}
Let $\la e,f\ra \cong U_{\mathrm{3A}}$ be a 3A-subalgebra of a VOA $V$.
\\
\textup{(1)}~ $U_{\mathrm{3A}}$ has 3 Ising vectors $e$, $f$ and $\tau_e f=\tau_fe$.
\\
\textup{(2)}~ $\abs{\tau_e \tau_f}=3$ and 
$\aut U_{\mathrm{3A}}=\la \tau_e,\tau_f\ra \cong \mathfrak{S}_3$.
\\
\textup{(3)}~ The 3A-algebra $U_{\mathrm{3A}}$ has a full subalgebra isomorphic to 
$L(\sfr{4}{5},0) \tensor L(\sfr{6}{7},0)$ where the central charges $c_3=4/5$ and $c_4=6/7$ are 
in the unitary series \eqref{eq:4.4}.
\end{thm}

\section{3-transposition groups and VOAs}

In this section we review some constructions of 3-transposition groups as 
automorphisms of vertex operator algebras.

\subsection{3-transposition groups and Ising vectors}

Let $E_V$ be the set of Ising vectors of a VOA $V$ of OZ-type and 
$E_{V}(\sigma)$ the subset of $E_V$ consisting of Ising vectors of $\sigma$-type on $V$.

\begin{prop}[\cite{M1}]\label{prop:5.1}
  If $e$, $f\in E_V(\sigma)$ then the type of the dihedral subalgebra $\la e,f\ra$ 
  is either 1A, 2A or 2B.
\end{prop}

Let $e$, $f\in E_V(\sigma)$.
Then their product in the Griess algebra looks as follows (cf.~\cite{M1,Ma}).
\begin{equation}\label{eq:5.1}
  ef=
  \begin{cases} 
    2e & \mbox{if $e=f$, i.e., $\la e,f\ra$ is of 1A-type,}
    \medskip\\
    0 & \mbox{if $(e\mymid f)=0$, i.e., $\la e,f\ra$ is of 2B-type,}
    \medskip\\
    \dfr{1}{4}(e+f-\sigma_e f) & \mbox{if $(e\mymid f)=2^{-8}$, i.e., $\la e,f\ra$ is of 2A-type.}
  \end{cases}
\end{equation}
Note that it follows from (4) of Theorem 4.8 that $\abs{\sigma_e\sigma_f}=3$ 
if $\la e,f\ra$ is of 2A-type, whereas $\abs{\sigma_e\sigma_f}\leq 2$ if $\la e,f\ra$ is 
of 2B-type (cf.~\cite{M1}).
We find that the relation above is the same as the local structure in \eqref{eq:3.1} 
of the Matsuo algebra with $\alpha=\beta=1/2$. 
Therefore, if we collect all the Ising vectors of $\sigma$-type, then we obtain 
a realization of a 3-transposition group as well as the associated Matsuo algebra 
in the Griess algebra.

\begin{thm}[\cite{M1,Ma}]\label{thm:5.2}
  Let $E_V(\sigma)$ be the set of Ising vectors of a VOA $V$ of OZ-type and 
  let $G_V=\la \sigma_e \in \aut V \mid e\in E_V(\sigma)\ra$.
  Then $G_V$ is a 3-transposition group.
  The Griess algebra of the subalgebra $\la E_V(\sigma)\ra$ of $V$ is isomorphic to 
  the non-degenerate quotient of the Matsuo algebra $B_{1/2,1/2}(G_V)$.
\end{thm}

Historically, the Matsuo algebra was formulated based on Theorem \ref{thm:5.2}.
Matsuo has classified the 3-transposition groups realized by Ising vectors of $\sigma$-type.

\begin{thm}[\cite{Ma, JLY}]
  Let $E_V(\sigma)$ be the set of Ising vectors of a VOA $V$ of OZ-type.
  The 3-transposition group $G_V=\la \sigma_e \in \aut V \mid e\in E_V(\sigma)\ra$ 
  is of symplectic type and a non-trivial connected component of $G_V$ is isomorphic 
  to one of the following:
  \[
    \mathfrak{S}_{n\geq 3},~~
    F^{\leq 2}{:}\mathfrak{S}_{n\geq 4},~~
    \mathrm{O}^+_{8\,\mathrm{or}\,10}(2),~~
    2^8{:}\mathrm{O}^+_8(2),~~
    \mathrm{O}^-_{6\,\mathrm{or}\,8}(2),~~
    2^6{:}\mathrm{O}^-_6(2),~~
    \mathrm{Sp}_{6\,\mathrm{or}\,8}(2),~~
    2^{6}{:}\mathrm{Sp}_6(2),
  \]
  where $F\cong 2^{2n}$ is the natural module for 
  $\mathrm{Sp}_{2n}(2)>\mathfrak{S}_{2n+2~\mathrm{or}~2n+1}$.
  Moreover, if $V$ is generated by $E_V(\sigma)$ as a VOA and $V$ is simple then 
  $V$ is uniquely determined by its Griess algebra, the non-degenerate quotient of 
  the Matsuo algebra $B_{1/2,1/2}(G_V)$ associated with $G_V$.
\end{thm}

In \cite{Ma}, the 3-transposition groups realizable in this manner were classified but 
the uniqueness of structures of VOAs realizing each 3-transposition group in Matsuo's list 
was open, and this uniqueness was established in \cite{JLY} under the assumption that 
the VOA is simple.
It turns out that all Matsuo algebras appearing in Theorem \ref{thm:5.2} are isomorphic to 
some subalgebras of the Griess algebras of $V_{\sqrt{2}R}^+$ associated with the root lattice 
$R$ of ADE-type.
Recall that a 3-transposition group $(G,I)$ is of symplectic type if and only if 
there is no $I$-subgroup isomorphic to $H\cong 3^{1+2}{:}2$ or its quotient $H/Z(H)\cong 3^2{:}2$
(cf.~Eq.~\eqref{eq:2.1}).
In Matsuo's classification, it is crucial to prove the absence of the $I$-subgroup 
isomorphic to $3^2{:}2$.
The proof of this part goes as follows.
Suppose we have a subalgebra $B$ isomorphic to $B_{1/2,1/2}(H)$ in the Griess algebra of $V$.
Then one can find a $c=c_2$ Virasoro sub VOA $L(\sfr{7}{10},0)$ and its module isomorphic to 
$L(\sfr{7}{10},\sfr{7}{10})$ in the sub VOA $\la B\ra$ of $V$.
However, by the list of possible highest weights in \eqref{eq:4.4},  
there is no irreducible highest weight representation of highest weight $7/10$ with central charge
$c_2=7/10$.
Thus 3-transposition groups realizable by Ising vectors of $\sigma$-type are only 
of symplectic type.

\subsection{\texorpdfstring{$W_3$}{W3}-algebra and \texorpdfstring{$\Z_3$}{Z3}-symmetry}

We explain another example of Matsuo algebras arising in VOA theory.
Here we consider the third value $c_3=4/5$ of the unitary series \eqref{eq:4.4}.
The Virasoro VOA $L(\sfr{4}{5},0)$ can be extended to a larger VOA 
$W_3(\sfr{4}{5})=L(\sfr{4}{5},0)\oplus L(\sfr{4}{5},3)$ called the \emph{$W_3$-algebra} 
at $c=4/5$ (cf.~\cite{FaZa}).
It is shown in (loc.~cit.) (see also \cite{KMY,M2}) that $W_3(\sfr{4}{5})$ is rational 
and has 6 irreducible representations over $\C$ as follows.
\begin{equation}\label{eq:5.2}
\begin{array}{cccccc}
  &W_3(\sfr{4}{5})=L(\sfr{4}{5},0\oplus 3),~~
  &L(\sfr{4}{5},\sfr{2}{5}\oplus \sfr{7}{5}),~~
  &L(\sfr{4}{5},\sfr{2}{3})^\pm,~~
  &L(\sfr{4}{5},\sfr{1}{15})^\pm
  \medskip\\
  \xi ~:& 1 & 1 & \zeta^{\pm 1} & \zeta^{\pm 1} 
\end{array}
\end{equation}
Here $L(c,h \oplus h')$ denotes $L(c,h)\oplus L(c,h')$ and $L(\sfr{4}{5},h)^\pm$ denotes 
inequivalent irreducible $W_3(\sfr{4}{5})$-modules which are isomorphic to $L(\sfr{4}{5},h)$ 
as $L(\sfr{4}{5},0)$-modules.
Let $\xi$ be the assignment of the powers of the cubic root $\zeta=e^{2\pii/3}$ of unity as above. 
Then it is known that $\xi$ preserves the fusion rules (cf.~\cite{FaZa,M2}).
Therefore, fusion rules of $W_3(\sfr{4}{5})$ have a $\Z_3$-symmetry.

Let $V$ be a VOA of OZ-type over $\R$.
A $c=4/5$ Virasoro vector $x\in V$ is called \emph{extendable} in $V$ 
if $\la x\ra\cong L(\sfr{4}{5},0)$ can be extended to a larger sub VOA 
$W_3(\sfr{4}{5})=L(\sfr{4}{5},0\oplus 3)$ inside $V$ 
(such an extension is known to be unique if exists).
Based on the $\Z_3$-symmetry of the fusion rules of $W_3(\sfr{4}{5})$-modules, 
one can define an automorphism $\xi_x$ of $\C V=\C\tensor_\R V$ associated to 
an extendable $c=4/5$ Virasoro vector $x\in V$ in the same manner as Miyamoto involutions.

\begin{lem}[\cite{M1}]\label{lem:5.4}
  Let $x\in V$ be an extendable $c=4/5$ Virasoro vector and $\xi_x$ 
  the automorphism of the complexified VOA $\C V$ defined as in \eqref{eq:5.2}.
  \\
  \textup{(1)}~ The map $\xi_x$ keeps the real space $V$ invariant and define 
  an element of $\aut V$.
  \\
  \textup{(2)}~ The $c=4/5$ Virasoro vector $x$ is of $\sigma$-type on the fixed point 
  subalgebra $V^{\la \xi_x\ra}$.
\end{lem}

Let $F_V$ be the set of extendable $c=4/5$ Virasoro vectors of $V$ and 
$F_V(\sigma)$ the subset of $F_V$ consisting of those of $\sigma$-type.
It follows from \eqref{eq:4.4}, \eqref{eq:5.2} and (2) of Lemma \ref{lem:5.4} 
that $x\in F_V(\sigma)$ if and only if $\xi_x=1$ on $V$.
Based on $c=4/5$ Virasoro vectors, we obtain another family of Matsuo algebras with 
$\alpha=2/5$ and $\beta=4/5$ as follows.

\begin{thm}[\cite{LY1}]
  Let $V$ be a simple VOA of OZ-type over $\C$ and let 
  $G_V=\la \sigma_x \in \aut V \mid x\in F_V(\sigma)\ra$.
  Then $G_V$ is a 3-transposition group.
  The Griess algebra of the subalgebra $\la F_V(\sigma) \ra$ of $V$ is isomorphic to 
  the non-degenerate quotient of the Matsuo algebra $B_{2/5,4/5}(G_V)$ associated with $G_V$.
\end{thm}

The classification of 3-transposition groups based on extendable $c=4/5$ Virasoro vectors 
of $\sigma$-type is open.
Some systematic examples are constructed in \cite{LY1} and it turns out that 
in contrary to the case of Ising vectors of $\sigma$-type, there exists an example of a VOA 
realizing a 3-transposition group of orthogonal type based on $c=4/5$ Virasoro vectors.

\begin{exam}[\cite{LY1}]
Let $K_{12}$ be the Coxeter-Todd lattice of rank 12 and consider the lattice VOA $V_{K_{12}}$.
There exists an order three automorphism $\nu$ of $V_{K_{12}}$ such that $X=V_{K_{12}}^{\la \nu\ra}$ 
is of OZ-type. 
Let $F_X(\sigma)$ be the set of extendable $c=4/5$ Virasoro vectors of $X$ which are of 
$\sigma$-type on $X$ and take $G_X=\la \sigma_x \mid x\in F_X(\sigma)\ra$. 
Then $G_X\cong {}^+\Omega_8^-(3)$.
It is worthy to mention that $K_{12}$ can be defined over the Eisenstein ring and 
its automorphism group preserving this complex structure is $3.{}^+\Omega_6^-(3)$ (cf.~\cite{CS}) 
so that there is a canonical subgroup $3^6{:}{}^+\Omega_6^-(3)$ in $\aut X$.
The interesting point here is that if we consider the automorphism group of 
the VOA $X=V_{K_{12}}^{\la \nu\ra}$ then the canonical symmetry $3^6{:}{}^+\Omega_6^-(3)$ can 
be extended to a larger symmetry ${}^+\Omega_8^-(3)$ which is highly non-trivial.
This is due to the fact that $X=V_{K_{12}}^{\la \nu\ra}$ has some extra $c=4/5$ extendable 
Virasoro vectors of $\sigma$-type.
For details, see Section 5.5 of \cite{LY1}.
\end{exam}

\begin{question}
  Is there an example of a VOA $V$ such that 
  $G_V=\la \sigma_x \in \aut V \mid x\in F_V(\sigma)\ra$ is of unitary type?
  If exists, then $B_{2/5,4/5}(2^{1+6}{:}H)$ appears as a Griess subalgebra 
  of such a VOA.
\end{question}

\subsection{3-transposition groups not arising from Matsuo algebras}

As we have seen in Theorem \ref{thm:5.2}, the 3-transposition groups generated by 
$\sigma$-involutions associated with Ising vectors of $\sigma$-type are always of symplectic type 
so that we cannot obtain 3-transposition groups other than of symplectic type as long as 
we consider Ising vectors of $\sigma$-type.
On the other hand, by the 6-transposition property in Sakuma's theorem, 
we cannot directly obtain 3-transposition groups if we use $\tau$-involutions associated 
with Ising vectors without any restriction.
However, if we collect suitable subsets of Ising vectors, then we can construct 
3-transposition groups generated by $\tau$-involutions.
Here we explain a general recipe of a construction of 3-transposition groups from 
Ising vectors which are not arising from Matsuo algebras.

Let $V$ be a VOA of OZ-type over $\R$ and 
let $E_V$ be the set of Ising vectors of $V$.
Suppose we have a pair $a$, $b\in E_V$ such that $(a\,|\,b)=13\cd 2^{-10}$ which 
is equivalent to that $\la a,b\ra$ is a 3A-algebra by (2) of Lemma \ref{lem:4.7}. 
We call such a pair a \emph{3A-pair}.
Fix a 3A-pair $a$, $b\in E_V$ and define
\begin{equation}\label{eq:5.3}
  I_{a,b}:=\{\, x\in E_V \mid (a\,|\,x)=(b\,|\,x)=2^{-5}\,\} .
\end{equation}
By Lemma \ref{lem:4.7},  $x\in I_{a,b}$ if and only if 
$\la a,x\ra\cong \la b,x\ra$ are 2A-algebras.
In this case one can also show that $\la a,b,x\ra=\la \sigma_x a, b\ra$ is a 6A-algebra 
(cf.~\cite{LY2}).

\begin{thm}[\cite{LY2}]\label{thm:5.8}
  Let $a$, $b\in E_V$ be a 3A-pair and define $I_{a,b}\subset E_V$ as above.
  \\
  \textup{(1)}~ If $x$, $y \in I_{a,b}$ then the type of $\la x,y\ra$ is either 1A, 2A or 3A. 
  In particular, one has $\abs{\tau_x\tau_y}\leq 3$ on $V$.
  \\
  \textup{(2)}~Set $G_V:=\la\, \tau_x \mid x\in I_{a,b} \,\ra$.
  Then $G_V$ is a 3-transposition group acting on the commutant subalgebra 
  $\com_V \la a,b\ra \subset V$ of $\la a,b\ra$ in $V$ whereas $G_V$ acts on $\la a,b\ra$ 
  trivially.
\end{thm}

Hereafter, $\com_V A$ denotes the commutant subalgebra of $A$ in $V$ (cf.~\cite{FZ}).
It follows from (2) of Theorem \ref{thm:4.8} that $G_V$ in Theorem \ref{thm:5.8} centralizes
$\la \tau_a,\tau_b\ra\cong \mathfrak{S}_3$.

\begin{exam}\label{exam:5.9}
  Let $V^\natural$ be the moonshine VOA \cite{FLM}.
  It is known that all the 3A-pairs in $E_{V^\natural}$ are mutually conjugate under the Monster 
  $\M=\aut V^\natural$ (cf.~\cite{C,M1}).
  Let $a$, $b\in V^\natural$ be Ising vectors such that $\la a,b\ra$ is a 3A-subalgebra 
  and take $I_{a,b}\subset E_{V^\natural}$ as in \eqref{eq:5.3}.
  Then we obtain $G_{V^\natural}=\la \tau_x \mid x\in I_{a,b}\ra\cong \Fi_{23}$ (cf.~\cite{LY2}).
  Indeed, Theorem \ref{thm:5.8} is hinted by the fact that 
  $\mathfrak{S}_3\times \Fi_{23}<\M$ where $\mathfrak{S}_3=\la \tau_a,\tau_b\ra$.
\end{exam}

As in the example above, we see that 3-transposition groups obtained by Theorem \ref{thm:5.8} 
include those of sporadic type.
However, this example strongly relies upon the Monster and there is no conceptual understanding 
of what kind of 3-transposition groups we can obtain from Theorem \ref{thm:5.8}.
Apart from this extreme example, we can simply say that the 3-transposition groups obtained by 
Theorem \ref{thm:5.8} are not arising from Matsuo algebras with $\beta=1/2$ as follows.
If $x$, $y\in I_{a,b}$ satisfy $\abs{\tau_x\tau_y}=3$ then $x$ and $y$ is also a 3A-pair.
In this case, the Griess algebra of the dihedral subalgebra $\la x,y\ra$ is 4-dimensional 
by Theorem \ref{thm:4.6} whereas subalgebras generated by two axes in Matsuo algebras are 
at most 3-dimensional. 
Therefore, 3-transposition groups obtained by Theorem \ref{thm:5.8} are beyond 
those arising from Matsuo algebras with $\beta=1/2$.

\subsection{3-transposition groups as homomorphic images}\label{sec:5.4}

We have seen a construction of 3-transposition group  based 
on $\tau$-involutions associated with Ising vectors in Theorem \ref{thm:5.8}.
Here we give another slightly different construction of 3-transposition groups.

Let $V$ be a VOA of OZ-type over $\R$ and $E_V$ the set of Ising vectors of $V$.
We consider all the 3A-pairs in $V$.
Set 
\begin{equation}\label{eq:5.4}
  T_V:=\{\, (x,y)\in {E_V}^2 \mid (x\,|\,y)=13\cd 2^{-10}\,\}
\end{equation}
and define
\begin{equation}\label{eq:5.5}
  u_{x,y}:=\dfr{448}{135}(x+y+\tau_x y)-\dfr{512}{405}
  (x+y+\tau_x y)_{(1)}(x+y+\tau_x y) \in \la x,y\ra \cong U_{\mathrm{3A}}.
\end{equation}
The vector $u_{x,y}$ above provides the $c=4/5$ Virasoro vector described 
in (3) of Theorem \ref{thm:4.9}.

\begin{lem}
If $(x,y)\in T_V$ then $u_{x,y}$ is an extendable $c=4/5$ Virasoro vector of $\la x,y\ra$.
\end{lem}

Fix a 3A-pair $(a,b)\in T_V$ and set
\begin{equation}\label{eq:5.6}
  J_{a,b}:=\l\{\, x \in V \mid \exists (x,y)\in T_V \mbox{ s.t. } u_{x,y}=u_{a,b}\,\r\},
  ~~
  G:=\la\, \tau_x \mid x\in J_{a,b}\,\ra .
\end{equation}
By Theorem \ref{thm:4.6}, $G$ itself is not a 3-transposition group but a 6-transposition group 
in general but we can extract a 3-transposition group from $G$ as follows.

\begin{thm}[\cite{HLY2}]\label{thm:5.11}
For a 3A-pair $(a,b)\in T_V$, define $J_{a,b}$ and $G$ as in \eqref{eq:5.6}.
Then the following restriction map is a group homomorphism.
\[
  \begin{array}{cccc}
    \psi_{a,b} :& G & \longrightarrow & \aut \com_V\la u_{a,b}\ra
    \\
    & g & \longmapsto & g|_{\com_V\la u_{a,b}\ra} 
  \end{array}
\]
The homomorphic image 
$\psi_{a,b}(G)=\la \,\psi_{a,b}(\tau_x) \mid x\in J_{a,b}\,\ra$ forms a 3-transposition group 
acting on the commutant subalgebra $\com_V \la u_{a,b}\ra$ of the Virasoro sub VOA 
$\la u_{a,b}\ra\cong L(\sfr{4}{5},0)$ of $V$.
\end{thm}

\begin{exam}\label{exam:5.12}
  If we apply Theorem \ref{thm:5.11} to $V^\natural$ then we obtain $G\cong 3.\Fi_{24}$ and 
  $\psi_{a,b}(G)\cong \Fi_{24}$ for any $(a,b)\in T_{V^\natural}$ since $\M=\aut V^\natural$ acts 
  on $T_{V^\natural}$ transitively (cf.~\cite{C,M1}).
  Note that $3.\Fi_{24}<\M$ but $\Fi_{24}\not<\M$ and $3.\Fi_{24}$ is a 6-transposition group but 
  not a 3-transposition group so that we need to take a homomorphic image
  to obtain a 3-transposition group in Theorem \ref{thm:5.11}.
\end{exam}

\subsection{4-transposition groups as homomorphic images}

In this subsection we explain a construction of 4-transposition groups similar to 
the one in the previous subsection.

Let $V$ be a VOA of OZ-type over $\R$ and $E_V$ the set of Ising vectors of $V$ as before.
Fix $e \in E_V$ and set
\begin{equation}\label{eq:5.7}
  I_e:=\{\, x\in E_V \mid (e\,|\,x)=2^{-5} \,\},~~~
  G:=\la \,\tau_x \in \aut(V) \mid x\in I_e\,\ra .
\end{equation}
By Theorem \ref{thm:4.6}, $G$ itself is a 6-transposition group in general.
Similar to Theorem \ref{thm:5.11}, we can extract a 4-transposition group from $G$ as follows.

\begin{thm}[\cite{HLY1}]\label{thm:5.13}
Take $e\in E_V$ and define $I_e$ and $G$ as in \eqref{eq:5.7}.
Then the following restriction map is a group homomorphism.
\[
\begin{array}{ccccc}
  \varphi_e &:& G &\longrightarrow& \aut \com_V \la e\ra 
  \\ 
  && g & \longmapsto & g|_{\com_V \la e\ra}
\end{array} 
\]
The homomorphic image $\varphi_e(G)=\la\, \varphi_e (\tau_x) \mid x\in I_e\,\ra$ 
is a 4-transposition group acting on the commutant subalgebra $\com_V \la e \ra$ of 
the Virasoro sub VOA $\la e\ra\cong L(\shf,0)$ of $V$.
\end{thm}

\begin{exam}\label{exam:5.16}
  If we apply Theorem \ref{thm:5.13} to $V^\natural$ then we obtain $G\cong 2.\mathbb{B}$ and 
  $\varphi_e (G)\cong \mathbb{B}$ for any $e\in E_{V^\natural}$ since $\M=\aut V^\natural$ acts 
  on $E_{V^\natural}$ transitively(cf.~\cite{C,M1}).
  Note that $2.\mathbb{B}<\M$ but $\mathbb{B}\not<\M$ and $2.\mathbb{B}$ is a 
  6-transposition group but not a 4-transposition group so that we need to take a homomorphic 
  image to obtain a 4-transposition group in Theorem \ref{thm:5.13}.
\end{exam}

\section{Open questions}

We list a few open questions about 3-transposition groups, Matsuo algebras and 
vertex operator algebras.

\begin{prob}
  Classify 3-transposition groups arising from Matsuo algebra $B_{2/5,4/5}(G)$ 
  based on $W_3(\sfr{4}{5})$. 
  Are they all of orthogonal type?
  To determine this, one needs to analyze subalgebras of the Matsuo algebra 
  $B_{2/5,4/5}(2^{1+6}{:}H)$ and determine whether there exists a subalgebra realizable 
  by a VOA or not.
\end{prob}

\begin{prob}
  Other than $(\alpha,\beta)=(1/2,1/2)$ and $(2/5,4/5)$, it is possible to realize 
  Matsuo algebras with $(\alpha,\beta)=(1/2,1/16)$ as Griess algebras of VOAs \cite{Ma,CH1,CH2}.
  It is shown in \cite{CH1} that the Matsuo algebra $B_{1/2,1/16}(H)$ for $H$ in \eqref{eq:2.1} 
  has a realization by a VOA so that 3-transposition groups of non-symplectic type 
  are possible to be realized as automorphism groups of VOAs via Matsuo algebras with 
  $(\alpha,\beta)=(1/2,1/16)$. 
  Some examples are constructed in \cite{CH2}.
  Give more examples of VOAs having Matsuo algebras with $(\alpha,\beta)=(1/2,1/16)$ 
  as Griess algebras and classify the 3-transposition groups realizable in this way if possible.
\end{prob}

\begin{prob}
  One may consider Majorana algebras \cite{Iv} and axial algebras \cite{HRS} 
  as axiomatizations of Griess algebras generated by Ising vectors.
  Can we generalize the construction of 3-transposition groups in Theorem \ref{thm:5.8} to
  arbitrary axial algebras?
  Also, can we generalize Theorems \ref{thm:5.11} and \ref{thm:5.13} to general axial algebras?
\end{prob}

\begin{prob}
  It seems possible to generalize Matsuo algebras to ``non-simply laced'' cases.
  In the original definition, all idempotents have the same spectrum and squared norm 
  but we can mix two or more different idempotents.
  For instance, when we consider the 3-transposition groups obtained by Theorems \ref{thm:5.8} 
  and \ref{thm:5.11} we observe some commutative algebras spanned by two kinds of idempotents 
  corresponding to $c_1=1/2$ and $c_4=6/7$ Virasoro vectors. 
  We do not include explicit examples here but we refer the interested readers to Section 4 
  of \cite{LY2}.
\end{prob}

\small

\end{document}